\def\obs#1{{\bf (*** #1 ***)} }

\def\obs#1{}     % Remova esta linha para rodar a versao 1
\documentclass[twoside,letterpaper,draft,11pt]{amsart}
\usepackage{amsmath,amsthm,latexsym,xspace,amscd,amssymb,mathtools,color,enumerate}

\usepackage[mathscr]{eucal}
\usepackage{amssymb}
\usepackage[all]{xy}
\makeatletter\renewcommand\theenumi{\@roman\c@enumi}\makeatother

\newtheorem{teo1}{Theorem}[section]

\newtheorem{lem1}[teo1]{Lemma}
\newtheorem{cor1}[teo1]{Corollary}
\newtheorem{prop1}[teo1]{Proposition}

\newtheorem{exe}[teo1]{Example}
\newtheorem{remark}[teo1]{Remark}

\newcommand{\m}{{}^{-1}}

\newcommand{\vu}{\vspace{.1cm}}
\newcommand{\vd}{\vspace{.2cm}}

\newcommand{\nod}{\noindent}

\newcommand{\id}{{\rm id}}

\newcommand{\A}{{A}}

\newcommand{\af}{\alpha}
\newcommand{\bt}{\beta}

\newcommand{\ta}{\tau}

\newcommand{\aft}{\af^\ast}

\newcommand{\G}{\mathcal{G}}

\def\ndv{\ {\mid \kern -0.7 em {\scriptstyle \not}} \ \ }

\def\nd{\ {\mid \kern -0.4 em {\scriptstyle \not}} \ \ }

\begin{document}

\thispagestyle{empty}

\title[On partial skew groupoids rings]{On partial skew groupoids rings}

\author[D. Bagio]{Dirceu Bagio}
\address{ Departamento de Matem\'atica, Universidade Federal de Santa Maria, 97105-900\\
	Santa Maria-RS, Brasil}
\email{bagio@smail.ufsm.br}

\author[A. Paques]{Antonio Paques }
\address{Instituto de Matem\'atica e Estat\'istica, Universidade federal de Porto Alegre, 91509-900\\
	Porto Alegre-RS, Brazil}
\email{paques@mat.ufrgs.br}

\author[H. Pinedo]{H\'ector Pinedo}
\address{Escuela de Matematicas, Universidad Industrial de Santander, Cra. 27 Calle 9 ´ UIS
	Edificio 45, Bucaramanga, Colombia}
\email{hpinedot@uis.edu.co}


\date{\today}
\begin{abstract} Given a partial action $\alpha$ of a connected groupoid $\G$ on an associative ring $A$ we investigate under what conditions the partial skew groupoid ring $A\star_{\alpha}\G$ can be realized as a partial skew group ring. In such a case applications concerning to the separability, semisimplicity and Frobenius property  of the ring extension $A\subset A\star_{\alpha}\G$ as well as to the artinianity of $\A\star_{\af}\G$ are given.
\end{abstract}

\maketitle

\setcounter{tocdepth}{1}

%\tableofcontents

\section{Introduction}

In this work we will consider partial actions of groupoids on rings. We are interested in studying the structure of the corresponding partial skew groupoids rings.

\vspace{.05cm}

Partial groupoids actions on rings were introduced in \cite{BP} and they are a natural generalization of partial group actions. It is well known that every groupoid is a direct sum of its connected component. A partial action of a groupoid on a ring $A$ is completely determined by the partial actions of its connected components on $A$. Thence, we can reduce the study of groupoid partial action on rings to the context of connected groupoids. 

\vspace{.05cm}

The structure of a connected groupoid is also well known. If $\G$ is a connected groupoid then $\G\simeq \G_0^2\times \G(x)$, where $\G_0^2$ is the coarse groupoid associated to the set $\G_0$ of the objects of $\G$ and $\G(x)$ is the isotropy group of an object $x$ of $\G$.

 \vspace{.05cm}

For a partial action $\alpha$ of a connected groupoid $\G$ on a ring $A$ we can construct the partial skew groupoid ring $A\star_\af \G$. If $\G_0$ is finite and $\alpha$ is unital then $A\star_\af \G$ is an associative and unital ring which is an extension of $A$.  \vspace{.05cm}

The partial skew groupoid rings have an important role in the partial Galois theory for groupoids as it is explicit in Theorem 5.3 of \cite{BP}. They also are examples of Leavitt path algebras, which are important in the theory of $C^{\ast}$-algebras (see Theorem 3.11 of \cite{GY}). In the last years, algebraic properties to the extension $A\subset A\star_\af \G$ have been studied. For example, the separability and semisimplicity properties of the extension $A\subset A\star_\af \G$ were studied in \cite{BPi} whereas in \cite{NOP} the authors investigate chain conditions between $A$ and $A\star_\af \G$.  

\vspace{.05cm}

Our purpose in this work is to study the following problem. Let $\G$ be a connected groupoid  such that $\G_0$ is finite and $\alpha$ a unital partial action of $\G$ on a ring $A$. Does the factorization $\G\simeq \G_0^2\times \G(x)$ induce a factorization of $A\star_\af \G$? Theorem 4.1 provides sufficient conditions for the answer to this question to be affirmative. Precisely, when $\alpha$ is a group-type partial action, we construct a  groupoid action $\beta$ of $\G_0^2$ on $A$ and a partial group action $\gamma$ of $\G(x)$ on $A\star_{\beta} \G_0^2$ and we prove that  $A\star_\af \G\simeq (A\star_{\beta} \G_0^2)\star_{\gamma}\G(x)$.

\vspace{.05cm}

We organize our work as follows. The background about groupoids is presented in Section 2. The topics of partial groupoid actions that will be used  are in Section 3.
In Section 4, we construct the actions $\beta$ and $\gamma$ which allow us to prove the factorization of $\A\star_{\af}\G$  mentioned in the previous paragraph. Applications of this result concerning to the separability, semisimplicity and Frobenius property of the extension $A\subset \A\star_{\af}\G$ as well as to the artinianity of $\A\star_{\af}\G$ are given in Section 5.

\subsection*{Conventions}\label{subsec:conv}
Throughout this work,  by ring we mean an associative and not necessarily unital ring. The center of a ring $A$ will be denoted by $C(A)$. We will denote the cardinality of a finite set $X$ by $|X|$.

\section{Groupoids}

We recall that a {\it groupoid} is a small category in which every morphism is an isomorphism. The set of the objects of a groupoid $\G$ will be denoted by $\G_0$.
If $g:x\to y$ is a morphism of $\G$ then $s(g)=x$ and $t(g)=y$ are called the {\it source} and the {\it target} of $g$ respectively. We will identify any object $x$ of $\G$ with its identity morphism, that is, $x=\id_x$.  The {\it isotropy group} associated to an object $x$ of $\G$ is the group $\G(x)=\{g\in \G:\,s(g)=t(g)=x\}$.\vspace{.05cm}

The composition of morphisms of a groupoid $\G$ will be denoted via concatenation. Hence, for $g,h\in \G$, there exists $gh$ if and only if $t(h)=s(g)$.
Notice that, if $g\in \G$ then $s(g)=g^{-1}g$ and $t(g)=gg\m$. Also, $s(gh)=s(h)$ and $t(gh)=t(g)$ for all $g,h\in \G$ with $t(h)=s(g)$.\vspace{.05cm}

A groupoid $\G$ is said to be {\it connected} if for any $x,y\in \G_0$ there exists a morphism $g\in \G$ such that $s(g)=x$ and $t(g)=y$, that is, the morphism $g$ connects the objects $x$ and $y$. It is well-known that any groupoid is a disjoint union of connected subgroupoids. In order to justify this fact, we consider the following equivalence relation on $\G_0$: for any $x,y\in\G_0$, $x\sim y$  if and only if there exists $ g\in\G$  such that $s(g)=x$ and $t(g)=y$.
Every equivalence class $X\in\G_0/\!\!\sim$ determines a full connected subgroupoid $\G_X$ of $\G$. The set of objects of $\G_X$ is $X$. The set ${\G_X}(x,y)$ of morphisms of $\G_X$ from $x$ to $y$ is equal to ${\G}(x,y)$, for all $x,y\in X$. By construction, $\G$ is the disjoint union of the subgroupoids $\G_X$, i.~e.
\begin{equation}\label{direct-sum}
\G=\dot\cup_{X\in \G_0/\!\sim}\G_X.
\end{equation}

For the convenience of the reader, we will prove a well-known result about the structure of connected groupoids. In order to do this, we need to introduce some extra notation.
Let $X$ be a nonempty set and $X^2=X\times X$. Then $X^2$ is a groupoid. The source and target maps of $X^2$ are, respectively, $s(x,y)=x$ and $t(x,y)=y$, for all $x,y\in X$. The rule of composition is given by: $(y,z)(x,y)=(x,z)$, for all $x,y,z\in X$. The groupoid $X^2$ is called the {\it coarse groupoid associated to $X$}.

\begin{prop1}\label{group:connec}
	Let $\G$ be a connected groupoid. Then $\G\simeq \G_0^2\times \G(x)$ as groupoids.
\end{prop1}
\begin{proof}
Let $x\in \G_0$ a fixed object of $\G$. For each $y\in \G_0$, we choose a morphism $\tau_y:x\to y$ of $\G$. We also choose $\tau_x=x$. Define $\varphi:\G\to \G_0^2\times\G(x)$
by $\varphi(g)= ((s(g),t(g)), g_x)$, where $g_x=\ta^{-1}_{t(g)}g\ta_{s(g)}$, for all $g\in \G$. It is straightforward to prove that $\varphi$ in a groupoid morphism.
Suppose that $\varphi(g)$ is an identity of $\G_0^2\times\G(x)$. Then, $\varphi(g)=((y,y),x)$ for some $y\in\G_0$. Hence, $s(g)=t(g)=y$ and $x=g_x=\ta_y^{-1}g\ta_y$ which implies that $g=\ta_y\ta_y^{-1}=x$. This ensures that $\varphi$ is injective. Given an element $((y,z),h)\in \G_0^2\times \G(x)$, consider $g=\tau_zh\tau\m_y\in \G$. Notice that $g_x=h$ and whence $\varphi(g)=((y,z),h)$, that is,  $\varphi$ is surjective, so an isomorphism of groupoids. \end{proof}

\section{Partial actions}
In this section we recall the notion of partial actions of groupoids. Some properties related to partial actions, that will be used later, are presented. The definition of group-type partial groupoid actions, which has a central role for our purposes, will be introduced.
\subsection{Partial groupoid action}
We recall from \cite{BP} that a \emph{partial action} of a groupoid $\G$ on  a ring $\A$ is a pair
$\af=(\A_g,\af_g)_{g\in \G}$  such that

\begin{enumerate}
	\item [(i)] $\A_g$ is an ideal of $\A_{t(g)}$ and $\A_{t(g)}$ is an ideal of $\A$, for all $g\in \G$,
	\item [(ii)]  $\af_g:\A_{g\m}\to \A_g$ is an isomorphism of rings, for all $g\in \G$,
	\item [(iii)] $\alpha_x=\id_{A_x}$, for all $x\in \G_0$,
	\item [(iv)] $\af_g\af_h\le \af_{gh}$, for all $g,h\in \G$ such that $t(h)=s(g)$.
\end{enumerate}
The condition (iv) means that $\alpha_{gh}$ is an extension of $\af_g\af_h$. Since the domain of $\af_g\af_h$ is $\af_h^{-1}(\A_{g^{-1}}\cap\A_h)$, it follows that (iv) is equivalent to
\[\quad{\rm (v)}\quad
\af_h^{-1}(\A_{g^{-1}}\cap\A_h)\subset \A_{(gh)^{-1}}\,\text{ and }\,\af_{gh}(a)=\af_g\af_h(a), \, \text{ for all }\, a\in \af_h^{-1}(\A_{g^{-1}}\cap\A_h).\]

The partial action $\af$ is said to be {\it global} if $\af_g\af_h=\af_{gh}$, for all $g,h\in \G$ such that $t(h)=s(g)$. Also, $\af$ is called {\it unital} if each $A_g$ is a unital ring, i.~e., there exists a central element $1_g$ of $A$ such that $A_g=A1_g$, for all $g\in \G$.

Now we recall Lemma 1.1 of \cite{BP} which give us some useful properties of partial actions that will be used in what follows.

\begin{lem1}\label{lem:BP} Let $\af=(\A_g,\af_g)_{g\in \G}$ be a partial action of a groupoid $\G$ on a ring $A$. Then:
\begin{enumerate}
	\item [${\rm(i)}$] $\alpha$ is global if and only if $A_g=A_{t(g)}$, for all $g\in \G$;
	\item [${\rm(ii)}$] $\af_{g\m}=\af\m_g$, for all $g\in \G$;
	\item [${\rm(iii)}$] $\af_g(\A_{g\m}\cap \A_h)=\A_{g}\cap \A_{gh}$, for all $g,h\in \G$ such that $t(h)=s(g)$.
\end{enumerate}
\end{lem1}

\begin{remark} \label{obs-pag}{\rm  Let $\af=(\A_g,\af_g)_{g\in \G}$ be a partial action of a groupoid $\G$ on a ring $A$. Notice  that $\af$ induces by restriction a partial action $\af_{(x)}=(\A_h,\af_h)_{h\in\G(x)}$ of the isotropy group $\G(x)$ on the ring $\A_x$, for each $x\in\G_0$.}
\end{remark}

\begin{remark}{\rm Let $\G$ be a groupoid. Using the decomposition of $\G$ given in \eqref{direct-sum}, it is straightforward to check that partial actions of $\G$ on  a ring $\A$ induce by restriction partial actions of $\G_X$ on $\A$, for all $X\in \G_0/\!\!\sim$. Conversely,  partial actions of $\G$ on $A$ are uniquely determined by partial actions of $\G_X$, $X\in \G_0/\!\!\sim$, on $A$ . Hence, we can reduce the study of partial groupoid actions to the connected case.}
\end{remark}

\subsection{Group-type partial groupoid action}

Let $\G$ be a connected groupoid, $x\in \G_0$ and $\mathcal{S}_x=\{h\in \G\,:\,s(h)=x\}$. Consider the following equivalence relation on $\mathcal{S}_x$:
$$g\sim_x l\,\,\ \text{ if and only if }\,\, t(g)=t(l), \qquad g,l\in \mathcal{S}_x.$$
A transversal $\tau(x)=\{\tau_{y}:y\in \G_0\}$ for $\sim_x$ such that $\tau_x=x$ will be called a \emph{transversal for $x$}. Hence,
$\ta_y:x\to y$ is a chosen morphism of $\G$, for each $y\in\G_0$ and $\ta_x=x$.

A partial action $\af=(A_g,\af_g)_{g\in \G}$ of a connected groupoid $\G$ on $A$ will be called {\it  group-type} if there exist $x\in \G_0$ and a transversal $\tau(x)=\{\tau_{y}:y\in \G_0\}$ for $x$ such that
\begin{align}
\label{cond1} A_{\tau\m_y}=A_x \ \ \text{and} \ \ A_{\tau_y}=A_y, \ \ \text{ for all } \ y\in\G_0.
\end{align}

\begin{remark}{\rm
(i) Notice that the notion of group-type partial action not depend on the choice of object $x$. Indeed, for another object $z$ of $\G$, consider $\tilde{\tau}_y:=\tau_{y}\tau\m_{z}$, for all $y\in \G_0$.
Clearly, $\tilde{\tau}(z)=\{\tilde{\tau}_y:y\in \G_0\}$ is a transversal for $z$. From \eqref{cond1} follows that  $\af_{\tau_{y}}\af_{\tau\m_{z}}=\af_{\tau_{y}\tau\m_{z}}=\af_{\tilde{\tau}_y}$. Thus, $A_{(\tilde{\tau}_y)\m}=A_z$ and $A_{\tilde{\tau}_y}=A_y$, for all $y\in \G_0$.
\vspace{.07cm}

\noindent (ii) We use the term ``group-type partial actions" since by Theorem \ref{teo-decomp}, proved in the next section for this kind of partial actions, the corresponding partial skew groupoid ring is indeed a partial skew group ring.}
\end{remark}

By Lemma \ref{lem:BP} (i), any global groupoid action is group-type. The converse is not true as we can see in the next example.
\begin{exe}\label{57}
	{\rm Let $\G=\{g,h,l,m,l^{-1},m^{-1}\}$ be the groupoid with objects $\G_0=\{x,y\}$ and the following composition rules
		\[ g^2=x,\quad h^{2}=y,\quad lg=m=hl,\quad g\in\G(x),\quad h\in\G(y)\,\,\text{ and }\,\, l,m:x\to y. \]
		The diagram bellow illustrates the structure of $\G$:
		\[\xymatrix{& x\ar[r]^{l}  &y \ar[d]^{h}\\
			& x\ar[u]^{g} \ar[r]^{m} & y} \]
Consider $A=\mathbb{C}e_1\oplus\mathbb{C}e_2\oplus\mathbb{C}e_3\oplus\mathbb{C}e_4$, where $\mathbb{C}$ denotes the complex number field, $e_ie_j=\delta_{i,j}e_i$ and $e_1+\ldots+e_4=1$. We define the following partial action $\alpha=\big(A_p,\af_p\big)_{p\in\G}$ of $\G$ on $A$:
		\begin{align*}
		&A_x=\mathbb{C}e_1\oplus\mathbb{C}e_2=A_{l\m},& & A_y=\mathbb{C}e_3\oplus\mathbb{C}e_4=A_l, &\\[.2em]
		& A_g=\mathbb{C}e_1=A_{g\m}=A_{m\m},& & A_m=A_h=\mathbb{C}e_3=A_{h\m}, &\
		\end{align*}
		and
		\begin{align*}
		&\af_x=id_{A_x},\ \ \ \af_y= id_{A_y},\ \  \ \af_g:ae_1\mapsto \overline{a}e_1, \ \ \ \af_h:ae_3\mapsto \overline{a}e_3,\ \  \ \af_m:ae_1\mapsto \overline{a}e_3, \\[.3em]
		&\ \ \af_{m\m}: ae_3\mapsto\overline{a}e_1,\ \ \ \af_l:ae_1+be_2\mapsto ae_3+be_4, \ \ \ \af_{l\m}:ae_3+be_4\mapsto ae_1+be_2,
		\end{align*}
		where $\overline{a}$ denotes the complex conjugate of $a$, for all $a\in\mathbb{C}$. Notice that $\af$ is a group-type (not global) partial action. Indeed, to see this it is enough to take the transversal $\ta(x)=\{\ta_x=x, \ta_y=l\}$ for $x$.}
\end{exe}

\section { The partial skew groupoid ring}

In this section, we will assume that  $\G$ is a connected groupoid such that $\G_0$ is finite, $x\in\G_0$ is a fixed object of $\G$ and $\af=(A_g,\af_g)_{g\in \G}$ is a unital partial action of $\G$ on a ring $A$ with $A_g=A1_g$, where $1_g$ is a central idempotent of $A$, for all $g\in\G$. We will also assume that $\alpha$ is group-type  and $\tau(x)=\{\tau_{y}:y\in \G_0\}$ is a transversal for $x$ such that \eqref{cond1} is satisfied. \vspace{.05cm}

The {\it partial skew groupoid ring} $A\star_\af\G$ associated to $\alpha$ is the set of all formal sums $\sum_{g\in \G}a_g\delta_g$, where $a_g\in A_g$, with the usual addition and multiplication induced by the following rule
\[
(a_g\delta_g)(a_h\delta_h)=
\begin{cases}
a_g\alpha_g(a_h1_{g^{-1}})\delta_{gh} &\text{if $s(g)=t(h),$}\\
0  &\text{otherwise},
\end{cases}
\]
for all $g, h\in \G$, $a_g\in A_g$ and $a_h\in A_h$.
The partial skew groupoid ring $A\star_\af\G$ is an associative ring. Since by assumption $\G_0$ is finite, $A\star_\af\G$ is  unital with identity  $1_{A\star_\af\G}=\sum_{y\in\G_0}1_y\delta_y$ (see $\S\,3$ of \cite{BFP} for more details). \vspace{.05cm}

As it was mentioned in the introduction section of \cite{NOP}, the ring structure of $A\star_\alpha \G$ only depends on the choice of the ideals $A_y$, $y\in \G_0$. Hence, we can choose $A$ to be any ring having the ideals as above described. In this sense, we will assume for the rest of this paper that
\[A=\oplus_{y\in \G_0} A_y.\]

\subsection{The main theorem}

In this subsection we will prove that the factorization of $\G$, given by Proposition \ref{group:connec}, induces a factorization of $A\star_\af\G$. Particularly, we obtain that  $A\star_\af\G$ is a partial skew group ring. In order to prove this result we will use some lemmas that will be proved in the sequel.

\begin{lem1}\label{xx}
The pair $\bt=(B_u,\bt_u)_{u\in\G_0^2}$, where  $B_u=A_{t(u)}$ and $\bt_u=\af_{\ta_{t(u)}}\af_{\ta\m_{s(u)}}$, is a global action of $\G_0^2$ on $A$.
\end{lem1}

\begin{proof}
	For any identity $e=(y,y)$ of $\G_0^2$ we have that $B_e=A_y$ and $\bt_e=\af_{\ta_{y}}\af_{\ta\m_{y}}$ is the identity map of $A_{\ta_y}=A_y$. Also, if $u=(y,z)$ and $v=(r,y)$ are elements in $\G_0^2$  then $uv=(r,z)$ and  $ \bt_u \bt_v=\af_{\ta_z}\af_{\ta\m_y}\af_{\ta_y}\af_{\ta\m_r}=\af_{\ta_z}\af_{\ta\m_r}= \bt_{uv}$.
\end{proof}

Thanks to Lemma \ref{xx} we can  consider the skew groupoid ring $C:=A\star_ {\bt}\G_0^2$.
In the sequel we will see that the group $\G(x)$ acts partially on $C$. Let $z\in \G_0$. Since $\af$ is group-type, it follows from \eqref{cond1} that $A_{\ta\m_z}=A_x$. Then, for all $h\in \G(x)$, $A_h\subset A_x$ and
\begin{align}\label{Czh}
C_{z,h}:=\af_{\ta_z}(A_h),
\end{align}
is well-defined. Hence, we can set
\begin{align}\label{Ch}
C_h:=\oplus_{u\in\G_0^2}C_{t(u),h}\delta_u.
\end{align}

\begin{lem1}\label{63} $C_h$ is a unital ideal of $C$, for all $h\in G(x)$. Moreover, $C_x=C$.
\end{lem1}
\begin{proof}
	Note that
	\begin{align*}
	C &=\oplus_{u\in \G_0^2}B_{u}\delta_u\\
	  &=\oplus_{u\in \G_0^2}A_{{t(u)}}\delta_u\\
	  &\overset{\mathclap{\eqref{cond1}}}{=} \oplus_{u\in \G_0^2}A_{\ta_{t(u)}}\delta_u \\
      &= \oplus_{u\in \G_0^2}\af_{\ta_{t(u)}}(A_{\ta\m_{t(u)}})\delta_u \\
      &\overset{\mathclap{\eqref{cond1}}}{=} \oplus_{u\in \G_0^2}\af_{\ta_{t(u)}}(A_x)\delta_u \\
      &\overset{\mathclap{\eqref{Czh}}}{=} \oplus_{u\in \G_0^2}C_{t(u),x}\delta_u \\
      &\overset{\mathclap{\eqref{Ch}}}{=}C_x.
      \end{align*}
Observe also that $1'_h=\sum_{z\in\G_0}\af_{\ta_z}(1_h)\delta_{(z,z)}$ is the identity element of $C_h$, for all $h\in\G(x)$. Indeed, let $u=(y,w)\in \G_0^2$ and $a\in C_{t(u),h}\delta_u=C_{w,h}\delta_u$. By \eqref{Czh}, there exists $a_h\in A_h$ such that $a=\af_{\ta_w}(a_h)\delta_{(y,w)}$ and consequently
\begin{align*}
a1'_h &=\sum_{z\in\G_0}\af_{\ta_w}(a_h)\delta_{(y,w)}\af_{\ta_z}(1_h)\delta_{(z,z)}\\
&=\af_{\ta_w}(a_h)\delta_{(y,w)}\af_{\ta_y}(1_h)\delta_{(y,y)}\\
&=\af_{\ta_w}(a_h) \bt_{(y,w)}(\af_{\ta_y}(1_h))\delta_{(y,w)}\\
&=\af_{\ta_w}(a_h)\af_{\ta_w}\af_{\ta\m_y}\af_{\ta_y}(1_h)\delta_{(y,w)} \qquad   \text{(see Lemma \ref{xx}))}  \\
&=\af_{\ta_w}(a_h)\af_{\ta_w}(1_h)\delta_{(y,w)}\\
&=\af_{\ta_w}(a_h)\delta_{(y,w)}\\
&=a.
\end{align*}
Similarly, $1'_ha=a$. Hence $1'_h$ is a central idempotent of $C$.  A straightforward calculation shows  that $C_h=1'_hC$. Thus, $C_h$ is an ideal of $C$.
\end{proof}

Let $(z,h)\in \G_0\times \G(x)$. We define $\gamma_{z,h}:C_{z,h\m}\to C_{z,h}$, $\af_{\ta_z}(a)\mapsto \af_{\ta_z}(\af_h(a))$, for all $a\in A_{h\m}$. Clearly $ \theta_{z,h}$ is a bijection. Moreover, these maps induce the following bijective map
\[ \gamma_h:C_{h\m}\to C_h,\quad  \gamma_h(\af_{\ta_{t(u)}}(a)\delta_u)= \gamma_{t(u),h}(a)\delta_u, \, \text{ for all } a\in A_{h\m} \text{ and } u\in \G_0^2.\]

\begin{lem1}\label{64}
The pair	$\gamma=(C_h, \gamma_h)_{h\in \G(x)}$ is a unital partial action of $\G(x)$ on $C$.
\end{lem1}	

\begin{proof}
By definition, $\gamma_x$ is the identity map of $C_x=C$. Note that $\gamma_h$ preserves the operation of multiplication. In fact, for all $a,b\in A_{h\m}$,
\begin{align*}
\hspace*{-1.5cm} \gamma_h((\af_{\ta_z}(a)\delta_{(y,z)})(\af_{\ta_y}(b)\delta_{(w,y)}))&= \gamma_h(\af_{\ta_z}(a) \aft_{(y,z)}(\af_{\ta_y}(b))\delta_{(w,z)})\\
&= \gamma_h(\af_{\ta_z}(a)\af_{\ta_z}(\af_{\ta\m_y}(\af_{\ta_y}(b)))\delta_{(w,z)})\\
&= \gamma_h(\af_{\ta_z}(ab)\delta_{(w,z)})\\
&=\af_{\ta_z}(\af_h(ab))\delta_{(w,z)}\\
&=\af_{\ta_z}(\af_h(a))\af_{\ta_z}(\af_h(b))\delta_{(w,z)}.
\end{align*}
On the other hand,
\begin{align*}
\hspace*{-1.2cm} \gamma_h(\af_{\ta_z}(a)\delta_{(y,z)}) \gamma_h(\af_{\ta_y}(b)\delta_{(w,y)})&=\af_{\ta_z}(\af_h(a))\delta_{(y,z)}\af_{\ta_y}(\af_h(b))\delta_{(w,y)}\\
&=\af_{\ta_z}(\af_h(a))\af_{\ta_z}(\af_{\ta\m_y}(\af_{\ta_y}(\af_h(b))))\delta_{(w,z)}\\
&=\af_{\ta_z}(\af_h(a))\af_{\ta_z}(\af_h(b))\delta_{(w,z)}.
\end{align*}
Hence, $\gamma_h$ is a ring isomorphism. It remains to show that $\gamma$ satisfies the condition (v) given in $\S\,3.1$.

Firstly, note that $\gamma_{l\m}(C_l\cap C_{h\m})\subset C_{(hl)\m}$, for all $h,l\in \G(x)$. Indeed, by definition, $\gamma_{l\m}$ is additive and whence
\begin{align*}
\hspace*{-1cm} \gamma_{l\m}(C_l\cap C_{h\m})&=\oplus_{u\in \G_0^2} \gamma_{l\m}(\af_{\ta_{t(u)}}(A_l\cap A_{h\m})\delta_u)\\
&=\oplus_{u\in \G_0^2}\af_{\ta_{t(u)}}(\af_{l\m}(A_l\cap A_{h\m}))\delta_u\\
&\subset \oplus_{u\in \G_0^2}\af_{\ta_{t(u)}}(A_{l\m}\cap A_{l\m h\m})\delta_u\\
&=C_{(hl)\m},
\end{align*}
where the last inclusion above holds because $\af_{(x)}$ is a partial action of $\G(x)$ on $A_x$ as defined in Remark \ref{obs-pag}.
Finally, let $c=\af_{\ta_z}(\af_{l\m}(a))\delta_{(y,z)}\in  \gamma_{l\m}(C_l\cap C_{h\m})$. Then
\begin{align*}
 \gamma_h( \gamma_l(c))&= \gamma_h( \gamma_l(\af_{\ta_z}(\af_{l\m}(a))\delta_{(y,z)}))\\
&= \gamma_h(\af_{\ta_z}(a)\delta_{(y,z)})\\
&=\af_{\ta_z}(\af_h(a))\delta_{(y,z)}.
\end{align*}
On the other hand, since $\af_{hl}=\af_h\af_l$ in $\af_{l\m}(A_l\cap A_{h\m})$ we have
\[ \gamma_{hl}(c)=\af_{\ta_z}(\af_{hl}(\af_{l\m}(a)))\delta_{(y,z)}=\af_{\ta_z}(\af_h(a))\delta_{(y,z)},\]
and consequently  $\gamma$ satisfies (v) of $\S\,3.1$. \end{proof}

By Lemma \ref{64} we can consider the partial skew group ring $(A\star_{\bt}\G_0^2)\star_ \gamma\G(x)$ and thus present the main result of this section which give us a factorization of the ring $A\star_\af\G$.

\begin{teo1} \label{teo-decomp}	
	$\,\,A\star_\af\G\simeq (A\star_{\bt}\G_0^2)\star_ \gamma\G(x).$
\end{teo1}

\begin{proof}
Consider the map $\varphi:A\star_\af\G\to (A\star_{\bt}\G_0^2)\star_ \gamma\G(x)$ given by
$ a\delta_g\mapsto a\delta_{(s(g),t(g))}\delta_{g_x}$,  where $g_x=\tau\m_{t(g)}g\tau_{s(g)}$.
In order to prove that $\varphi$ is a ring isomorphism we proceed by a series of steps.

\vu
	
\noindent	\emph{Step 1: $\varphi$ is well defined.}

\vu
	
\noindent By Lemma \ref{xx}, $A_g\subseteq A_{{t(g)}}=B_{(s(g),t(g))}$, for all $g\in\G$. Hence,  we only need to show that $a\delta_{(s(g),t(g))}\in C_{g_x}$, for all $a\in A_g$.
Notice that
	\begin{align*}
	\af_{\ta_{t(g)}}(A_{g_x})&=\af_{\ta_{t(g)}}(A_{\ta\m_{t(g)}g\ta_{s(g)}})\\
	&=\af_{\ta_{t(g)}}(A_{\ta\m_{t(g)}g\ta_{s(g)}}\cap A_x)\\
	&=A_{g\ta_{s(g)}}\cap A_{\ta_{t(g)}x} \qquad (\text{by Lemma \ref{lem:BP} {\rm (iii)}})\\
	&=A_{g\ta_{s(g)}}\cap A_{\ta_{t(g)}}\\
    &\overset{\mathclap{\eqref{cond1}}}{=}A_{g\ta_{s(g)}}\cap A_{t(g)}\\
	&=A_{g\ta_{s(g)}} \qquad \qquad\qquad (A_{g\ta_{s(g)}}\subseteq A_{t(g\ta_{s(g)})}=A_{t(g)} ).
	\end{align*}
Since $A_{g\m}=A_{g\m}\cap A_{s(g)}\overset{\eqref{cond1}}{=} A_{g\m}\cap A_{\ta_{s(g)}}$ we have $A_g=\af_g(A_{g\m}\cap A_{\ta_{s(g)}})=A_g\cap A_{g\ta_{s(g)}}$. Hence $a\in A_g\subseteq A_{g\ta_{s(g)}}=\af_{\ta_{t(g)}}(A_{g_x})$ which implies $a\delta_{(s(g),t(g))}\in C_{g_x}$ by \eqref{Czh} and \eqref{Ch}.

\vd

\noindent \emph{Step 2: $\varphi$ is a ring homomorphism.}

\vu	

\noindent It is enough to prove that $\varphi$ preserves the operation of multiplication. Let $g,h\in\G$ such that $s(g)=t(h)$. It is easy to see that $(gh)_x=g_xh_x$. Hence
\begin{align*}
\varphi((a\delta_g)(b\delta_h))&=\varphi(a\af_g(b1_{g\m})\delta_{gh})\\
&=a\af_g(b1_{g\m})\delta_{(s(gh),t(gh))}\delta_{(gh)_x}\\
&=a\af_g(b1_{g\m})\delta_{(s(h),t(g))}\delta_{g_xh_x},
\end{align*}
for all $a\in A_g$ and $b\in A_h$. On the other hand
\begin{align*}
\varphi(a\delta_g)\varphi(b\delta_h)&=(a\delta_{(s(g),t(g))}\delta_{g_x})(b\delta_{(s(h),t(h))}\delta_{h_x})\\
&=(a\delta_{(s(g),t(g))})(\gamma_{g_x}(b\delta_{(s(h),t(h))}1'_{g\m_x}))\delta_{g_xh_x}.
\end{align*}
As in Step 1, we have that $A_h\subseteq A_{h\ta_{s(h)}}=\af_{\ta_{t(h)}}(A_{h_x})$. Since $b\in A_h$, there is $b'\in A_{h_x}$ such that $b=\af_{\ta_{t(h)}}(b')$ and
\begin{align*}
b\delta_{(s(h),t(h))}1'_{g\m_x}&=\af_{\ta_{t(h)}}(b')\delta_{(s(h),t(h))}\sum_{z\in\G_0}\af_{\ta_z}(1_{g\m_x})\delta_{(z,z)}\\
&=\af_{\ta_{t(h)}}(b')\delta_{(s(h),t(h))}\af_{\ta_{s(h)}}(1_{g\m_x})\delta_{(s(h),s(h))}\\
&=\af_{\ta_{t(h)}}(b')\af_{\ta_{t(h)}}(1_{g\m_x})\delta_{(s(h),t(h))}\\
&=\af_{\ta_{t(h)}}(b'1_{g\m_x})\delta_{(s(h),t(h))}.
\end{align*}
Hence,
\begin{align*}
\varphi(a\delta_g)\varphi(b\delta_h)&=(a\delta_{(s(g),t(g))})(\gamma_{g_x}(b\delta_{(s(h),t(h))}1'_{g\m_x}))\delta_{g_xh_x}\\
&=(a\delta_{(s(g),t(g))})(\gamma_{g_x}(\af_{\ta_{t(h)}}(b'1_{g\m_x})\delta_{(s(h),t(h))}))\delta_{g_xh_x}\\
&=(a\delta_{(s(g),t(g))})(\af_{\ta_{t(h)}}(\af_{g_x}(b'1_{g\m_x}))\delta_{(s(h),t(h))})\delta_{g_xh_x}\\
&=a\af_{\ta_{t(g)}}\af_{\ta\m_{s(g)}}\af_{\ta_{t(h)}}\af_{g_x}(b'1_{g\m_x})\delta_{(s(h),t(g))}\delta_{g_xh_x}\\
&=a\af_{\ta_{t(g)}}\af_{g_x}(b'1_{g\m_x})\delta_{(s(h),t(g))}\delta_{g_xh_x}\\
&\overset{\mathclap{\eqref{cond1}}}{=}a\af_{\ta_{t(g)}}(\af_{g_x}(b'1_{g\m_x})1_{\ta\m_{t(g)}})\delta_{(s(h),t(g))}\delta_{g_xh_x}\\
&=a\af_{\ta_{t(g)}g_x}(b'1_{(\ta_{t(g)}g_x)\m})1_{\ta_{t(g)}}\delta_{(s(h),t(g))}\delta_{g_xh_x}\qquad (\text{by (v) of $\S\,$3.1})\\
&=a\af_{g\ta_{s(g)}}(b'1_{(g\ta_{s(g)})\m})1_{\ta_{t(g)}}\delta_{(s(h),t(g))}\delta_{g_xh_x}\\
&\overset{\mathclap{\eqref{cond1}}}{=}a\af_{g\ta_{s(g)}}(b'1_{(g\ta_{s(g)})\m})1_g\delta_{(s(h),t(g))}\delta_{g_xh_x}\\
&=a\af_g(\af_{\ta_{s(g)}}(b'1_{\ta\m_{s(g)}})1_{g\m})\delta_{(s(h),t(g))}\delta_{g_xh_x}\qquad (\text{by (v) of $\S\,$3.1})\\
&\overset{\mathclap{\eqref{cond1}}}{=}a\af_g(\af_{\ta_{t(h)}}(b')1_{g\m})\delta_{(s(h),t(g))}\delta_{g_xh_x}\\
&=a\af_g(b1_{g\m})\delta_{(s(h),t(g))}\delta_{g_xh_x}\\
&=\varphi((a\delta_g)(b\delta_h)).
\end{align*}

\noindent \emph{Step 3: $\varphi$ is injective.}
\vu

\noindent Let $v=\sum_{g\in\G}a_g\delta_g\in A\star_\af\G$ such that $\varphi(v)=0$. Then
$$0=\sum_{g\in\G}a_g\delta_{(s(g),t(g))}\delta_{g_x}=\sum_{h\in\G(x)}\sum_{\substack{g\in \G\\ g_x=h}} a_g\delta_{(s(g),t(g))}\delta_h$$
Since $C\star_ \theta\G(x)$ is a direct sum, it follows that
\begin{align}\label{injective}
\sum_{\substack{g\in \G\\ g_x=h}}a_g\delta_{(s(g),t(g))}=0,\quad \text{for all } h\in \G(x).
\end{align}
Consider $h\in\G(x)$ and $g,g'\in\G$ such that $g_x=g'_x=h$. It is straightforward to check that $(s(g),t(g))=(s(g'),t(g'))$ if and only if $g=g'$. Therefore \eqref{injective} holds if and only if $a_g=0$, for all $g\in \G$. Thus $v=0$.

\vd

\noindent \emph{Step 4: $\varphi$ is surjective.}

\vu

It is enough to check that given any element of the type $\af_{\ta_z}(a)\delta_{(y,z)}\delta_h$, with $h\in\G(x)$ and $a\in A_h$, there exists an element $w\in A\star_\af\G$ such that $\varphi(w)=\af_{\ta_z}(a)\delta_{(y,z)}\delta_h$. To do that observe that
\begin{align*}
\af_{\ta_z}(a)\in\af_{\ta_z}(A_h)&=\af_{\ta_z}(A_h\cap A_x)\\
&\overset{\mathclap{\eqref{cond1}}}{=}\af_{\ta_z}(A_h\cap A_{\ta\m_z})\\
&=A_{\ta_zh}\cap A_{\ta_z}\qquad\quad (\text{by Lemma \ref{lem:BP} (iii)})\\
&\overset{\mathclap{\eqref{cond1}}}{=}A_{\ta_zh}\cap A_z\\
&=A_{\ta_zh}\qquad\quad  (\text{because}\,\ A_{\ta_zh}\subseteq A_{t(\ta_zh)}=A_{t(\ta_z)}=A_z).
\end{align*}
Therefore, for  $g=\ta_zh\ta\m_y$ we have $t(g)=t(\ta_z)=z$, $s(g)=s(\ta\m_y)=y$ and  $$\af_{\ta_z}(a)\in A_{\ta_zh}= A_{\ta_{t(g)}h}\overset{(\ast)}{=}\af_{\ta_{t(g)}}(A_h)\overset{(\ast\ast)}{\subseteq}\af_{\ta_{t(g)}}(A_{h\ta\m_{s(g)}})\overset{(\ast)}{=} A_{\ta_{t(g)}h\ta\m_{s(g)}}=A_{\ta_zh\ta\m_y}=A_g,$$ where $(\ast)$ is ensured by
$$\af_{\ta_{t(g)}}(A_h)=\af_{\ta_{t(g)}}(A_h\cap A_x)\overset{\eqref{cond1}}{=}\af_{\ta_{t(g)}}(A_h\cap A_{\ta\m_{s(g)}})=A_{\ta_{t(g)}h},$$ and $(\ast\ast)$ by
$$A_h=\af_h(A_{h\m}\cap A_x)\overset{\eqref{cond1}}{=}\af_h(A_{h\m}\cap A_{\ta\m_{s(g)}})=A_h\cap A_{h\ta\m_{s(g)}}\subseteq A_{h\ta\m_{s(g)}}.$$

\nod Now, taking $w=\af_{\ta_z}(a)\delta_g$ we are done.	
\end{proof}

\begin{remark} \label{ob-global} {\rm Since global actions are group-type actions, the factorization given in Theorem \ref{teo-decomp} holds for all unital global action $\alpha$ of $\G$ on $A$. In such a case, the partial group action $\gamma$ of $\G(x)$ on $A_x$ is indeed a global action and consequently $A\star_{\af}\G$ is a skew group ring.}
	
\end{remark}

\section{Applications}

The aim of this section is to present some applications of Theorem \ref{teo-decomp}. In what follows, $\G$ is connected and $\G_0$ is finite. The partial action $\alpha$, the ring $A$ and the transversal $\tau(x)$ are assumed as in the previous section. Also, $\beta$ is the global action of $\G_0^2$ on $A$ given in Lemma \ref{xx} and $\gamma$ is the partial action of $\G(x)$ on $A\star_{\bt}\G_0^2$ given in Lemma \ref{64}. \smallbreak

Note that $\varphi\colon \A\to \A\star_{\af}\G$, $\ a\mapsto \sum_{y\in \G_0}(a1_y)\delta_y$, is a monomorphism of rings and whence $\A\star_{\af}\G$ is a ring extension of $A$.
By Theorem \ref{teo-decomp}, $A\subset \A\star_{\beta}\G_0^2\subset (\A\star_{\beta}\G_0^2)\star_{\gamma}\G(x)\simeq \A\star_{\af}\G$. Therefore, we will investigate some properties of the extension $A\subset \A\star_{\af}\G$ using the intermediate extensions and the results known for partial group actions.

\subsection{Separability} In this subsection we will study the separability property to the ring extension $A\subset A\star_{\af}\G$.
We recall that a  unital ring extension $R\subset S$ is called {\it separable} if the multiplication map $m:S\otimes_R S\to S$ is a splitting epimorphism of $S$-bimodules. This is equivalent to saying that there exists an element $x\in S\otimes_RS$ such that $sx=xs$, for all $s\in S$, and $m(x)=1_S$. Such a element $x$ is usually called {\it an idempotent of separability} of $S$ over $R$. \vspace{.1cm}

Throughout this subsection, we will assume that $\G$ is {\it finite}. As in \cite{BPi}, consider the maps $t_{y,z}:A\to A$ and $t_{z}:A\to A$ given  by
\[ t_{y,z}(a)=\sum_{g\in \G(y,z)} \af_g(a1_{g\m}), \quad t_z(a)=\sum_{y\in \G_0} t_{y,z}(a),\,\,\quad y,z\in \G_0,\,\,a\in A.\]
Particularly, if $\G$ is a group then $\G_0=\{x\}$ and $t_x:A\to A$ is the trace map for partial group actions as defined in Section 2 of \cite{DFP}.\vspace{.1cm}

Now, we  recall Theorem 4.1 of \cite{BPi} which will be useful for our purposes.
\begin{teo1} \label{teo:BPin}
	The ring extension $A\subset A\star_{\af}\G$ is separable if and only if  there is an element $a$ in the center $C(A)$ of $A$ such that $t_{z}(a)=1_z$, for all $z\in \G_0$.\qed
\end{teo1}

In order to apply Theorem \ref{teo-decomp} to determine when the extension $A\subset A\star_{\af}\G$ is separable, we consider the separability problem for the extensions $A\subset A\star_{\beta}\G_0^2$ and $A\star_{\beta}\G_0^2\subset (A\star_{\beta}\G_0^2)\star_{\gamma}\G(x)$.

\begin{lem1}\label{lem:sep1}
	The extension $A\subset A\star_\beta\G_0^2$ is separable if and only if there exists $a\in C(A)$ such that $\sum_{z\in\G_0}\af_{\ta\m_{z}}(a1_z)=1_x$.
\end{lem1}
\begin{proof}
Let $z\in \G_0$ and $a\in A$. Then
\begin{align*}
t_z(a)&=\sum_{y\in \G_0}t_{y,z}(a)=\sum_{y\in \G_0}\sum_{u\in \G_0^2(y,z)}\beta_u(a1_{u^{-1}})\\
      &=\sum_{y\in \G_0}\beta_{(y,z)}(a1_{y})=\sum_{y\in \G_0}\alpha_{\tau_z}\alpha_{\tau\m_y}(a1_{y})\\
      &=\alpha_{\tau_z}(\sum_{y\in \G_0}\alpha_{\tau\m_y}(a1_{y})).
\end{align*}
Consequently $t_z(a)=1_z$ if and only if $\sum_{y\in \G_0}\alpha_{\tau\m_y}(a1_{y})=1_x$ and the result follows by Theorem \ref{teo:BPin}. \end{proof}

\begin{lem1}\label{lem:centro}
	$C(A\star_\beta\G_0^2)=\left\{\sum_{z\in\G_0} \alpha_{\tau_z}(a_x)\delta_{(z,z)}\,:\,a_x\in C(A_x)\right\}.$
\end{lem1}
\begin{proof}
	Let $a_x\in C(A_x)$ and $\Lambda=\sum_{z\in\G_0} \alpha_{\tau_z}(a_x)\delta_{(z,z)}$. Given $(y,w)\in \G_0^2$ and $a_w\in A_w$ we have that
	\begin{align*}
	a_w\delta_{(y,w)}\cdot \Lambda=a_w\delta_{(y,w)}\cdot\alpha_{\tau_y}(a_x)\delta_{(y,y)}=a_w\alpha_{\tau_w}(a_x)\delta_{(y,w)}
	\end{align*}
	and
	\begin{align*}
	\Lambda\cdot a_w\delta_{(y,w)}=\alpha_{\tau_w}(a_x)\delta_{(w,w)}\cdot a_w\delta_{(y,w)}=\alpha_{\tau_w}(a_x)a_w\delta_{(y,w)}.
	\end{align*}
	Since $a_x\in C(A_x)$ and $\alpha_{\tau_w}$ is an isomorphism it is clear that $\alpha_{\tau_w}(a_x)\in C(A_w)$.  Thus, $\Lambda\in C(A\star_\beta\G_0^2)$. Conversely, consider $\Lambda=\sum_{y,z\in\G_0} a_{(y,z)}\delta_{(y,z)}\in C(A\star_\beta\G_0^2)$, with $a_{(y,z)}\in A_{(y,z)}=A_z$, for all $y\in \G_0$.
	From $\Lambda\cdot 1_w\delta_{(w,w)}= 1_w\delta_{(w,w)}\cdot \Lambda$, it follows that
	\[\sum_{z\in\G_0}a_{(w,z)}\delta_{(w,z)}=\sum_{y\in\G_0}a_{(y,w)}\delta_{(y,w)},\,\, \text{ for all } w\in \G_0.\]
	Hence $a_{(y,z)}=0$ if $y\neq z$ and whence $\Lambda=\sum_{z\in\G_0} a_{(z,z)}\delta_{(z,z)}$. Moreover, for all $y,w\in \G_0$ and $a\in A_w$,
	\[\Lambda\cdot a\delta_{(y,w)}=a_{(w,w)}a\delta_{(y,w)}\,\,\,\text{ and }\,\,\, a\delta_{(y,w)}\cdot \Lambda=a\alpha_{\tau_w}(\alpha_{\tau\m_y}(a_{(y,y)}))\delta_{(y,w)}.\]
	Therefore
	\begin{align}\label{equ:comu}
	a_{(w,w)}a=a\alpha_{\tau_w}(\alpha_{\tau\m_y}(a_{(y,y)})), \,\,\,\text{ for all }\,\, y,w\in \G_0,\,\,a\in A_w.
	\end{align}
	When $a=1_w$ we obtain that $a_{(w,w)}=\alpha_{\tau_w}(\alpha_{\tau\m_y}(a_{(y,y)}))$, for all $y,w\in \G_0$. Particularly, $a_{(w,w)}=\alpha_{\tau_w}(a_{(x,x)})$ for all $w\in \G_0$.
	Given  $b\in A_x$, consider $a:=\alpha_{\tau_w}(b)\in A_w$. By \eqref{equ:comu}, $\alpha_{\tau_w}(a_{(x,x)}b)=\alpha_{\tau_w}(ba_{(x,x)})$. Thus $a_{(x,x)}b=ba_{(x,x)}$ and consequently $a_{(x,x)}\in C(A_x)$.
\end{proof}

\begin{lem1}\label{lem:sep2} The extension  $A\star_\beta\G_0^2\subset (A\star_\beta\G_0^2)\star_{\gamma} \G(x)$ is separable if and only if the extension $A_x\subset A_x\star_{\alpha_{(x)}} \G(x)$ is separable.
\end{lem1}
\begin{proof}
Let $\Lambda\in C(A\star_\beta\G_0^2)$. By Lemma \ref{lem:centro}, $\Lambda=\sum_{z\in \G_0}\alpha_{\tau_z}(a_x)\delta_{(z,z)}$, with $a_x\in C(A_x)$.
Notice that
\begin{align*}
t_x(\Lambda)&=\sum_{h\in \G(x)}\gamma_h(\sum_{z\in \G_0}\alpha_{\tau_z}(a_x)\delta_{(z,z)}\cdot \sum_{w\in \G_0}\alpha_{\tau_w}(1_{h\m})\delta_{(w,w)} )\\
          &=\sum_{h\in \G(x)}\sum_{z\in \G_0}\gamma_h(\alpha_{\tau_z}(a_x)\delta_{(z,z)}\cdot \alpha_{\tau_z}(1_{h\m})\delta_{(z,z)} )\\
          &=\sum_{h\in \G(x)}\sum_{z\in \G_0}\gamma_h(\alpha_{\tau_z}(a_x1_{h\m})\delta_{(z,z)} )\\
          &=\sum_{h\in \G(x)}\sum_{z\in \G_0}\alpha_{\tau_z}(\alpha_h(a_x1_{h\m}))\delta_{(z,z)}\\
          &=\sum_{z\in \G_0}\alpha_{\tau_z}(\sum_{h\in \G(x)}\alpha_h(a_x1_{h\m}))\delta_{(z,z)}\\
          &=\sum_{z\in \G_0}\alpha_{\tau_z}(t_x(a_x))\delta_{(z,z)}.\\
\end{align*}
Hence, there is $\Lambda\in C(A\star_\beta\G_0^2)$ such that $t_x(\Lambda)=1=\sum_{z\in \G_0} 1_z\delta_{(z,z)}$ if and only if there is $a_x\in C(A_x)$ such that $t_x(a_x)=1_x$. it Then the result follows from Theorem \ref{teo:BPin} and Theorem 3.1 of \cite{BLP}.
\end{proof}

\begin{teo1}\label{teo:sep} The following statements hold:
	\begin{enumerate}
		\item [${\rm(i)}$] if $A_x\subset A_x\star_{\af_{(x)}}\G(x)$ is a separable extension and there exists $a\in C(A)$ such
 that $\sum_{z\in\G_0}\af_{\ta\m_{z}}(a1_z)=1_x$ then $A\subset A\star_{\af}\G$ is a separable extension; \smallbreak
 \item [${\rm(ii)}$] if $A\subset A\star_{\af}\G$ is a separable extension then $A_x\subset A_x\star_{\af_{(x)}}\G(x)$ also is.
\end{enumerate}
\end{teo1}
\begin{proof}
\noindent (i) By Lemmas \ref{lem:sep1} and \ref{lem:sep2}, $A\subset A\star_\beta\G_0^2$ and $A\star_\beta\G_0^2\subset (A\star_\beta\G_0^2)\star_{\gamma} \G(x)$ are separable extensions. By the transitivity of the separability property (see Proposition 2.5 of \cite{HS}), $A\subset (A\star_\beta\G_0^2)\star_{\gamma} \G(x)$ is separable. Thus the result follows from Theorem \ref{teo-decomp}. \vspace{.1cm}

\noindent (ii) From Theorem \ref{teo-decomp} and Proposition 2.5 of \cite{HS} it follows that $A\star_\beta\G_0^2\subset (A\star_\beta\G_0^2)\star_{\gamma} \G(x)$ is a separable extension. Hence, Lemma \ref{lem:sep2} implies that $A_x\subset A_x\star_{\af_{(x)}}\G(x)$ is separable.
\end{proof}

\begin{cor1} \label{cor-1}Suppose that $|\G_0|1_A$ is invertible in $A$ and $(|\G_0|1_A)^{-1}=n1_A$, with $n\in \mathbb{N}$. If $A_x\subset A_x\star_{\af_{(x)}}\G(x)$ is separable then $A\subset A\star_{\af}\G$ is separable.
\end{cor1}	
\begin{proof} Let $a=n1_A\in C(A)$. Then
\begin{align*}
\sum_{z\in\G_0}\af_{\ta\m_{z}}(a1_z)&=\sum_{z\in\G_0}\af_{\ta\m_{z}}(n1_z)=n\sum_{z\in\G_0}\af_{\ta\m_{z}}(1_z)\\
                                    &=n|\G_0|1_x=1_x,
                                    \end{align*}
and whence the result follows by Theorem \ref{teo:sep} (1).
\end{proof}

\begin{exe}
	{\rm Let $\G$, $A$, $\alpha$ and $\ta(x)$ be as in Example \ref{57}. Consider $a_x=\frac{1}{2}e_1+e_2\in A_x$. Observe that $t_x(a_x)=\af_x(a_x)+\af_g(a_x1_{g\m})=a_x+\frac{1}{2}e_1=e_1+e_2=1_x$. It follows from Theorem 3.1 of \cite{BLP} that $A_x\subset A_x\star_{\af_{(x)}} \G(x)$ is a separable extension. Moreover, if $a=\frac{1}{2}(e_1+e_2+e_3+e_4)\in A$ then
		\[\alpha_{\tau\m_x}(a1_x)+\alpha_{\tau\m_y}(a1_y)=\alpha_x(a1_x)+\alpha_{l\m}(a1_l)=e_1+e_2=1_x.\]
	By Theorem \ref{teo:sep} (1), the extension $A\subset A\star_{\alpha}\G$ is separable.	
	}
\end{exe}

\subsection{Semisimple extension}
In this subsection we investigate the semisimplicity for the ring extension $A\subset A\star_{\af}\G$. Recall from \cite{HS} that a ring extension $R\subset S$ is called left (right) {\it semisimple} if any left (right) $S$-submodule $N$ of a left (right) $S$-module $M$ having an $R$-complement
in $M$, has an $S$-complement in $M$.
\vspace{.15cm}

For the convenience of the reader, we recall Proposition 2.6 of \cite{HS}.

\begin{prop1} \label{sep-ss}
If $R\subset S$ is a separable ring extension then $R\subset S$ is a  left (right) semisimple ring extension. \qed
\end{prop1}

\begin{teo1}\label{teo-ss} If $\G$ is {\it finite} and
\begin{enumerate}
	\item [${\rm(i)}$] there exists $a_x\in C(A_x)$ such that $t_x(a_x)=1_x$ and \vspace{0.15cm}
	\item [${\rm(ii)}$] there exists $a\in C(A)$ such that $\sum_{z\in\G_0}\af_{\ta\m_{z}}(a1_z)=1_x$,
\end{enumerate}
then $A\subset A\star_{\af}\G$ is a left (right) semisimple extension.
\end{teo1}

\begin{proof} Using (i), we obtain from Theorem 3.1 of \cite{BLP} that $A_x\subset A_x\star_{\af_{(x)}}\G(x)$ is separable. Then, by Theorem \ref{teo:sep} (i) the extension
	 $A\subset A\star_{\af}\G$  is separable and so the result follows by Proposition \ref{sep-ss}.
\end{proof}

\subsection{Frobenius extension}
In this subsection we will prove that $A\subset A\star_{\af}\G$ is a Frobenius extension. We recall that a ring extension $R\subset S$ is called {\it Frobenius} if there exist an element $\Delta=\sum_{i=1}^{n}s_{i1}\otimes s_{i2}\in S\otimes_R S$ and an $R$-bimodule map $\varepsilon:S\to R$ such that $\Delta s=s\Delta$, for all $s\in S$, and $\sum_{i=1}^{n}\varepsilon(s_{i1})s_{i2}=\sum_{i=1}^{n}s_{i1}\varepsilon(s_{i2})=1$. More details on Frobenius extensions can be seen, for example, in \cite{CIM} or \cite{K}.  \vspace{.1cm}

Firstly, we note that the natural inclusion $A\subset A\star_{\beta} \G_0^2 $, given by $a\mapsto \sum_{z\in\G_0}(a1_z)\delta_{(z,z)}$, induces the following $(\A,\A)$-bimodule structure on $A\star_{\beta} \G_0^2$:
\begin{align*}\label{aabim}
&a\cdot (a_z\delta_{(y,z)})=aa_z\delta_{(y,z)},\quad\quad
(a_z\delta_{(y,z)})\cdot a=a_z\bt_{(y,z)}(a1_y)\delta_{(y,z)},&
\end{align*}
for all $(y,z)\in\G_0^2$, $a_z\in A_z$ and $a\in \A$.

\begin{teo1} \label{teo-frobenius} If $\G$ is finite then $A\subset A\star_{\af}\G$ is a Frobenius extension.
\end{teo1}
\begin{proof}
Let $\Delta=\sum_{z\in \G_0}1_z\delta_{(z,z)}\otimes 1_z\delta_{(z,z)}\in (A\star_{\beta}\G_0^2)\otimes_{A} (A\star_{\beta}\G_0^2)$ and
$\varepsilon:A\star_{\beta}\G_0^2\to A$	given by \[\varepsilon(a_z\delta_{(y,z)})=\left\{\begin{array}{lc}
a_z, 	& \text{ if } y=z\\
0,& \text{ otherwise.}
\end{array}
\right. \]
It is straightforward to check that $\varepsilon$ is an $(A,A)$-bimodule map. Also, note that
\[(a_w\delta_{(y,w)}) \Delta=a_w\delta_{(y,w)}=\Delta(a_w\delta_{(y,w)}),\]
for all $(y,w)\in \G_0^2$ and $a_w\in A_w$. Since
\[\sum_{z\in \G_0}\varepsilon(1_z\delta_{(z,z)})(1_z\delta_{(z,z)})=\sum_{z\in \G_0}1_z\delta_{(z,z)}\varepsilon(1_z\delta_{(z,z)})
=\sum_{z\in \G_0}1_z\delta_{(z,z)}=1_{A\star_{\beta}\G_0^2},\]
it follows that $A\subset A\star_{\beta}\G_0^2$ is a Frobenius extension. Notice that $\G(x)$ is finite because $\G$ is finite. Then, by Theorem 3.6 of \cite{BLP}, $A\star_{\beta}\G_0^2\subset (A\star_{\beta}\G_0^2)\star_{\gamma}\G(x)$ is a Frobenius extension. As Frobenius extension is a transitive notion (see e.~g. \cite[pg. 6]{K}), we obtain from Theorem \ref{teo-decomp} that $A\subset A\star_{\af}\G$ is Frobenius.
\end{proof}

\subsection{Artinianity} The artinian property for partial skew groupoids rings was studied in \cite{NOP}. In our context, using Theorem \ref{teo-decomp}, we obtain the following refinement of Theorem 1.3 of \cite{NOP}.

\begin{teo1} The partial skew groupoid ring $A\star_{\af}\G$ is artinian if and only if $A$ is artinian and $A_h=\{0\}$,  for all but finitely many $h\in \G(x)$. 
\end{teo1}
\begin{proof} Assume that $A\star_{\af}\G$ is artinian. By Theorem 1.3 of \cite{NOP}, $A$ is artinian and $A_g=\{0\}$ for all but finitely many $g\in \G$. Particularly, $A_h=\{0\}$ for all but finitely many $h\in \G(x)$. 
	
For the converse, consider $h\in \G(x)$. By $\eqref{Czh}$ and \eqref{Ch}, $C_h\neq \{0\}$ if and only if there is $z\in \G_0$ such that $\af_{\tau_z}(A_h)\neq \{0\}$. Consequently, $C_h=\{0\}$, for all but finitely many $h\in \G(x)$. Also, since $\G_0^2$ is finite and $\A$ is artinian it follows from Theorem 1.3 of \cite{NOP} that $\A\star_{\beta}\G_0^2$ is artinian. Using again the Theorem 1.3 of \cite{NOP}, we conclude that $(\A\star_{\beta}\G_0^2)\star_{\gamma}\G(x)$ is artinian and the result follows by Theorem \ref{teo-decomp}.
\end{proof}

\end{document}